\documentclass[10pt]{amsart}
\usepackage{amsmath,amscd}
\usepackage{amsbsy}
\usepackage{amssymb}
\usepackage{amscd,amsthm}
\usepackage[all,cmtip]{xy}
\usepackage[russian,english]{babel}
\usepackage[utf8]{inputenc}

\newtheorem{thm}{Theorem}\numberwithin{thm}{section}

\newtheorem{cor}[thm]{Corollary}

\newtheorem{prob}[thm]{Problem}

\newtheorem*{con2}{Conjecture}

\begin{document}
\begin{center}
\huge{On the diophantine equation $An!+Bm!=f(x,y)$}\\[1cm] %and a result of Erd\"os and Obl\'ath}\\[1cm]
\end{center}
\begin{center}

\large{Sa$\mathrm{\check{s}}$a Novakovi$\mathrm{\acute{c}}$}\\[0,5cm]
{\small September 2023}\\[0,5cm]
\end{center}
%\end{center}
%\begin{center}
%\begin{otherlanguage*}{russian}
%	\emph{ЗА APИАНУ.}	
%\end{otherlanguage*}
%\end{center}
{\small \textbf{Abstract}. 
Erd\"os and Obl\'ath proved that the equation $n!\pm m!=x^p$ has only finitely many integer solutions. More general, under the ABC-conjecture, Luca showed that $P(x)=An!+Bm!$ has finitely many integer solutions for polynomials of degree $\geq 3$. For certain polynomials of degree $\geq 2$, this result holds unconditionally. We consider irreducible homogeneous $f(x,y)\in \mathbb{Q}[x,y]$ of degree $\geq 2$ and show that there are only finitely many $n,m$ such that $An!+Bm!$ is represented by $f(x,y)$. As corollaries we get alternative proofs for the unconditional results of Luca. We also discuss the case of certain reducible $f(x,y)$. Furthermore, we study equations of the form $n!!m!!=f(x,y)$ and $n!!m!!=f(x)$.}

%we show that this implies even more, namely that a diophantine equation of the form $bA^n=P(x)$ with $P(x)$ a polynomial in $\mathbb{Z}[x]$ of degree $d\geq 2$ and $b$ and $A$ fixed has only finitely many integer solutions $(n,x)$. %In this way we generalize some previously known result.%We also give examples of w-helices by which we obtain solutions of Markov-type equations. 
\begin{center}
\tableofcontents
\end{center}
\section{Introduction}
Diophantine quations involving factorials have a long and rich history. For example Brocard \cite{BR} and independently Ramanujan \cite{RA} asked to find all integer solutions for $n!=x^2-1$. Up to now this is an open problem, known as Brocard's problem. It is believed that the equation has only the three solutions $(x,n)=(5,4), (11,5)$ and $(71,7)$. Overholt \cite{O} observed that a weak form of Szpiro's-conjecture implies that Brocard's equation has finitely many integer solutions. Using the ABC-conjecture Luca \cite{L} proved that diophantine equations of the form $n!=P(x)$ with $P(x)\in \mathbb{Z}[x]$ of degree $d\geq 2$ have only finitely many integer solutions with $n>0$. If $P(x)$ is irreducible, Berend and Harmse \cite{BH1} showed unconditionally that $P(x)=H_n$ has finitely many integer solutions where $H_n$ are highly divisible sequences which also include $n!$. Furthermore, they proved that the same is true for certain reducible polynomials. %In the present work, we consider diophantine equations of the form $bn_1!A^{n_1}\cdots n_r!A_r^{n_r}=P(x)$ with fixed non-zero integer $b$ and fixed positive integers $A_1,...,A_r$ and prove an analogous statement for such equations. %Note that Siegel \cite{S} proved that $bA^n+D=x^2$ has only finitely many integer solutions. 

Without assuming the ABC-conjecture, Berend and Osgood \cite{BO} showed that for arbitrary $P(x)$ of degree $\geq2$ the density of the set of positive integers $n$ for which there exists an integer $x$ such that $P(x)=n!$ is zero. %We believe that the arguments presented in \emph{loc.cit}. can also be applied in more general situations, implying, for instance, that the density of the set of positive integers solving $n_1!A^{n_1}\cdots n_r!A_r^{n_r}=P(x)$ is zero. 
Further progress in this direction was obtained by Bui, Pratt and Zaharescu \cite{BPZ} where the authors give an upper bound on integer solutions $n\leq N$ to $n!=P(x)$. %For a detailed overview on results about, for instance, Ramanujan--Nagell type equations $bA^n+D=x^2$ and exponential diophantine equations in general we refer to \cite{BMS}, \cite{SH} and \cite{ST}. 
Of course, there are several polynomials $P(x)$ for which $P(x)=n!$ is known to have either very few integer solutions or none (see for instance \cite{EO}, \cite{DAB} and \cite{PS}). %But finding all the integer solutions for $x^2-1=n!$ is still an unsolved problem, known as Brocard's problem \cite{BR}. 
Berndt and Galway \cite{BG} showed that the largest value of $n$ in the range $n<10^9$ for which Brocard's equation $x^2-1=n!$ has a positive integer solution is $n=7$. Matson \cite{MA} extended the range to $n<10^{12}$ and Epstein and Glickman \cite{EG} to $n<10^{15}$. %Overholt \cite{O} observed that a weak form of Szpiro's-conjecture actually implies that Brocard's equation has only finitely many integer solutions. 

Starting from Brocard's problem there are also studied variations or generalizations of $x^2-1=n!$ (see for instance \cite{DMU}, \cite{KL}, \cite{MU}, \cite{MUT} \cite{TY}). For instance Ulas \cite{MU} studied, among others, diophantine equations of the form $2^nn!=y^2$ and proved that the Hall conjecture (which is a special case of ABC-conjecture) implies that the equation has only finitely many integer solutions. Note that $2^nn!$ can also be formulated using the notation for the Bhargava factorial $n!_S$. Note that $2^nn!=n!_S$, with $S=\{2n+b| n\in \mathbb{Z}\}$. We do not recall the definition of the Bhargava factorial and refer to \cite{BH} or \cite{WT} instead. Other diophantine equations involving factorials have proved more tractable. For example, Erd\"os and Obl\'ath \cite{EO} proved that $n!\pm m!=x^p$ has only finitely many integer solutions. This result has been generalized to equations of the form $An!+Bm!=P(x)$ where $P(x)$ is a polynomial with rational coefficients of degree $\geq 2$. Luca \cite{L2} proved that if the degree is $\geq 3$, then the ABC-conjecture implies that $An!+Bm!=P(x)$ has finitely many integer solutions, except when $A+B=0$. In this case there are only finitely many solutions with $n\neq m$. For polynomials of the form $P(x)=a_dx^d+a_{d-3}x^{d-3}+\cdots +a_0$ with $a_d\neq 0$ and $d\geq 2$ Gawron \cite{MG} showed that there a finitely many integer solutions to $n!+m!=P(x)$, provided ABC holds.

In \cite{NO} the author started to study diophantine equations of the following type: let $g(x_1,...,x_r)\in\mathbb{Z}[x_1,...,x_r]$ and $f$ be polynomials where either $f\in\mathbb{Q}[x]$ or $f\in\mathbb{Q}[x,y]$. Consider the equation
\begin{center}
	$g(x_1,...,x_r)=f$,
\end{center}
where for any of the $x_i$ in $g$ we may also plug in $A^n$ or $n!_S$. Several results concerning the existence of finitely many integer solutions were proved by the author \cite{NO} for $g(x_1,...,x_r)=bx_1\cdots x_r$ and certain $f$. For details we refer to \emph{loc.cit.} In the above equation, we are mainly interested in $f\in\mathbb{Q}[x]$ or $f\in\mathbb{Q}[x,y]$. The reason for this is the following: it is known that a positive integer is represented by $f(x,y,z)=x^2+y^2+z^2$ if it is not of the form $4^l(8k+7)$, where $l,k$ are non-negative integers. So one can check for instance that there are infinitely many $n$ such that $n!$ or $A^nn!$ are not of the form $4^l(8k+7)$. So for certain $f\in \mathbb{Z}[x,y,z]$ we have an easy strategy to check whether there are infinitely many solutions. In four variables, it is known that $x^2+y^2+z^2+w^2$ represents any positive integer.

%We think that this is a reasonable generalization. Of course, formulated in this generality there are plenty of $g$ and $f$ such that the diophantibe equation has infinitely many (positive) integer solutions. 
Some interesting (exponential) diophantine equations are of the above form. Below we give only a few examples:
\begin{itemize}
	\item[(i)] \emph{superelliptic equations}.
	\item[(ii)] \emph{Erd\"os-Obl\'ath type equations}: take $g(x_1)=x_1$ and plug in $n!$ and let $g(x,y)=x^p\pm y^p$ or take $g(x_1,x_2)=x_1\pm x_2$ and plug in $n!$ respectively $m!$ and let $f(x)=x^p$.
	\item[(iii)] \emph{Thue-Mahler equation}: take $g(x_1,...,x_r)=x_1\cdot x_2\cdots x_r$ and let $x_i=p_i^{n_i}$.
	\item[(iv)] \emph{Thue-equation}: take $g(x_1,...,x_r)=m$ and let $f(x,y)$ be a homogeneous polynomial of degree $\geq 3$.
	\item[(v)] \emph{Brocard's problem}: take $g(x_1)=x_1$ and plug in $n!$ and let $f(x)=x^2-1$. More generally, let $f(x)$ be any polynomial of degree $\geq 2$ and we get the diophantine equation considered in \cite{L}.
	\item[(vi)] \emph{Ramanujan-Nagell type equation}: take $g(x_1)=bx_1+D$ and plug in $A^n$ for some positive fixed $A$ and let $f(x)=x^2$.
	\item[(vii)] \emph{Fermat equation}: take $g(x_1,x_2)=x_1^n+x_2^n$ and let $f(x)=x^n$.
	\item[(viii)] \emph{generalizations of Brocard's problem}, see \cite{WT}: take $g(x_1)=x_1$ and plug in $n!$ and let $f(x,y)$ be any irreducible binary form of degree $\geq 2$. 
	\item[(ix)] \emph{generalizations of Brocard's problem}, see \cite{DA}: take $g(x_1,x_2)=-x_1^2+x_2$ and plug in $n!$ and let $f(x,y)=x^2+y^2-A$.
\end{itemize}

In the examples (i) to (viii) from above, conditionally or unconditionally there are only finitly many integer solutions. In some situations the exact number of integer solutions is known. %So from a structural point of view one can ask for a geometric characterization of the polynomials $g$ and $f$ (or the variety defined by $g-f=0$) such that the diophantine equation $g=f$ has finitely many solutions. % In general, this is a deep problem, tackeld in diophantine geomety. 

In this note we want to study two types of euations: first, we want to study the case $g(x_1,x_2)=Ax_1+Bx_2$, where we plug in $n!$ and $m!$ respectively. Formulated in the notation of the Bhargava factorial, we plug in $n!_{\mathbb{Z}}$ respectively $m!_{\mathbb{Z}}$. %or $n!_S$, where $S=\{An+b | n\in\mathbb{Z}\}$ for some fixed positive integer $A$. %The main results of the present note are the following:
Therefore, we want to study diophantine equations of the form
\begin{center}
	$An!+Bm!=f(x)$
\end{center}
and 
\begin{center}
	$An!+ Bm!=f(x,y)$,
\end{center}
where $f(x)$ and $f(x,y)$ are polynomials with rational coefficients and $A,B$ non-zero fixed integers. If $n=m$ and $A+B=0$, then $f(x)=0$ respectively $f(x,y)=0$. Hence, we get infinitely many solutions. They all have $n=m$ and $f(x)=0$ respectively $f(x,y)=0$. If $n=m$ and $A+B\neq 0$, we have the diophantine equations 
\begin{center}
	$(A+B)n!=f(x)$ \textnormal{and} $(A+B)n!=f(x,y)$.
\end{center}
Now these equations have been studied in \cite{NO}. So we may assume $n>m$. The second type of equations we are interested in involve the double factorial $n!!$, which is closely related to a certain Bhargava factorial. Equations involving the double factorial have been studied by Ulas \cite{MU}. In the present note we also want to generalize some of the results in \cite{MU} by considering $g(x_1,...,x_r)=bx_1\cdots x_r$ and pluging in the double factorial. Here $b$ is a non-zero integer. More precisely, we study the equations:
\begin{center}
	$bn_1!!\cdots n_r!!=f(x)$ \textnormal{and} $bn_1!!\cdots n_r!!=f(x,y)$.
\end{center}
Note that if $n=2m$ we have $n!!=2^mm!=m!_S$ where $S=\{2m+b | m\in \mathbb{Z}\}$. And if $n=2m+1$ we have $n!!=\frac{(2m+1)!}{2^mm!}=\frac{(2m+1)!_{\mathbb{Z}}}{m!_S}$.\\

The proofs of Theorems 1.9 and 1.10 below need the ABC-conjecture. So we recall its content (see \cite{LA}). For a non-zero integer $a$, let $N(a)$ be the \emph{algebraic radical}, namely $N(a)=\prod_{p|a}{p}$. Note that 
\begin{eqnarray}
N(a)=\prod_{p|a}{p}\leq \prod_{p\leq a}{p}< 4^a,
\end{eqnarray}
where the last inequality follows from a Chebyshev-type result in elementary prime number theory and is called the Finsler inequality.
%\begin{con2}[Weak form of Szpiro's-conjecture]
%There exists some constant $s>0$ such that for mutually prime integers $A,B$ and $C$  with $A+B=C$ the inequality 
%\begin{eqnarray}
%|ABC|<N(ABC)^{s}
%\end{eqnarray}
%holds.
%\end{con2} 
%A generalization of the above conjecture is the following:
\begin{con2}[ABC-conjecture]
For any $\epsilon >0$ there is a constant $K(\epsilon)$ depending only on $\epsilon$ such that whenever $A,B$ and $C$  are three coprime and non-zero integers with $A+B=C$, then 
\begin{eqnarray}
\mathrm{max}\{|A|,|B|,|C|\}<K(\epsilon)N(ABC)^{1+\epsilon}
\end{eqnarray}
holds.
\end{con2} 
%The ABC-conjecture from above implies the weak form of Szpiro's-conjecture.\\
%\begin{prop}
%Let $P(x)\in\mathbb{Z}[x]$ a polynomial of degree $d\geq 2$ and $A$ a non-zero integer with $\mathrm{log}(A)>2d+1$. The ABC-conjecture implies that $A^n=P(x)$ has only finitely many integer solutions.
%\end{prop}
%We reprove Siegel's result for $bD^n+E=x^2$ by assuming that Szpiro's-conjecture holds true.
%\begin{prop}
%The weak form of Szpiro's-conjecture implies that $bD^n+E=x^2$ for fixed positive $b,D$ and $E$ has only finitely many integer solutions.
%\end{prop}

%In this note we want to study the case $g(x_1,x_2)=Ax_1+Bx_2$, where we plug in $n!$ and $m!$ respectively. Formulated in the notation of the Bhargava factorial, we plug in $n!_{\mathbb{Z}}$ respectively $m!_{\mathbb{Z}}$. %or $n!_S$, where $S=\{An+b | n\in\mathbb{Z}\}$ for some fixed positive integer $A$. %The main results of the present note are the following:
%Therefore, we want to study diophantine equations of the form
%\begin{center}
%	$An!+Bm!=f(x)$
%\end{center}
%and 
%\begin{center}
%	$An!+ Bm!=f(x,y)$,
%\end{center}
%where $f(x)$ and $f(x,y)$ are polynomials with integer coefficients and $A,B$ non-zero fixed integers. Furthermore, we want to study equations of the form
%\begin{center}
%	$n!!=f(x)$ \textnormal{and} $n!!=f(x,y)$.
%\end{center}
%Here $n!!$ denotes the double factorial. Equations involving the double factorial have been studied by Ulas \cite{MU}.
%For a better readibility and to keep it clearer, we prove the results in the present work for $q=r$. However, we would like to note that all results also hold for $q<r$.
We have the following results:
\begin{thm}Let $f(x,y)\in\mathbb{Q}[x,y]$ be an irreducibel binary form of degree $\geq 2$. Then there are only finitely many $n,m$ such that $An!+ Bm!$ is represented by $f(x,y)$, except when $A+B=0$. In this case there are only finitely many $n,m$ with $n\neq m$. If the degree is $\geq 3$, then the equation 
	\begin{eqnarray*}
		An!+Bm!=f(x,y)
	\end{eqnarray*}
	has only finitely many integer solutions, except when $A+B=0$. In this case there are only finitely many integer solutions with $n\neq m$.
\end{thm}
\begin{thm}
	Let $f(x,y)\in\mathbb{Q}[x,y]$ be a polynomial with factorization 
\begin{center}
	$f(x,y)=f_1(x,y)^{e_1}\cdots f_u(x,y)^{e_u}$,
\end{center}
where the  $f_i(x,y)$ are distinct irreducible homogeneous polynomials of degree $d_i$. 
If $\mathrm{min}\{d_1e_1,...,d_ue_u\}>1$,  
then there are only finitely many $n,m$ such that $An!+Bm!$ is represented by $f(x,y)$, except when $A+B=0$. In this case there are only finitely many $n,m$ with $n\neq m$. If the degree is $\geq 3$, then the equation 
\begin{eqnarray*}
	An!+Bm!=f(x,y)
\end{eqnarray*}
has only finitely many integer solutions, except when $A+B=0$. In this case there are only finitely many integer solutions with $n\neq m$.
\end{thm}
Theorems 1.1 and 1.2 generalize the main result of Luca \cite{L2} to certain binary forms. Consequently, we find Luca's unconditional results as corollaries:
\begin{cor}
	Let $f(x)\in\mathbb{Q}[x]$ be an irreducible polynomial of degree $\geq 2$. Then the diophantine equation $An!+Bm!=f(x)$ has only finitely many integer solutions, except when $A+B=0$. In this case there are only finitely many integer solutions with $n\neq m$.
\end{cor}
\begin{cor}
	Let $f(x)\in\mathbb{Q}[x]$ be a polynomial with factorization 
	\begin{center}
		$f(x)=f_1(x)^{e_1}\cdots f_u(x)^{e_u}$,
	\end{center}
	where $f_i(x)$ are distinct irreducible polynomials of degree $d_i$. If $\mathrm{min}\{d_1e_1,...,d_ue_u\}>1$, 
	then the equation $An!+ Bm!=f(x)$ has only finitely many integer solutions, except when $A+B=0$. In this case there are only finitely many integer solutions with $n\neq m$.
\end{cor}

\begin{thm}
	Let $f(x,y)\in\mathbb{Q}[x,y]$ be an irreducibel binary form of degree $d>r$ and $b$ a non-zero integer. Then there are only finitely many $(n_1,...n_r)$ such that $bn_1!!\cdots n_r!!$ is represented by $f(x,y)$. If the degree is $\geq 3$, then the equation 
\begin{eqnarray*}
bn_1!!\cdots n_r!!=f(x,y)
\end{eqnarray*}
has only finitely many integer solutions.
\end{thm}
\begin{thm}
	Let $b$ be a non-zero integer and $f(x,y)\in\mathbb{Q}[x,y]$ a polynomial with factorization 
	\begin{center}
		$f(x,y)=f_1(x,y)^{e_1}\cdots f_u(x,y)^{e_u}$,
	\end{center}
	where the  $f_i(x,y)$ are distinct irreducible homogeneous polynomials of degree $d_i$. Now if $d_i\geq2$ and $d_1e_1+\cdots +d_ue_u>r$ or if $\mathrm{min}\{d_1e_1,...,d_ue_u\}>r$, then there are only finitely many $(n_1,...n_r)$ such that $bn_1!!\cdots n_r!!$ is represented by $f(x,y)$. If the degree is $\geq 3$, then the equation 
	\begin{eqnarray*}
		bn_1!!\cdots n_r!!=f(x,y)
	\end{eqnarray*}
	has only finitely many integer solutions.
\end{thm}
Note that if $d\leq r$ one usually has infinitely many solutions. Indeed, let $f(x,y)=a_dx^d+a_{d-1}x^{d-1}y+\cdots + a_{d-1}xy^{1}+a_0y^d$ and assume $a:=a_d+\cdots +a_0\neq 0$. If $x=y$, we obtain the equation
\begin{center}
	$bn_1!!\cdots n_r!!=f(x,y)=ax^d$.
\end{center}
Now if all $n_i$ are even then \cite{NO}, Theorem 1.1 yields  that there are infinitely many solutions for $b\neq 0$. 
\begin{cor}
	Let $f(x)\in\mathbb{Q}[x]$ be an irreducible polynomial of degree $d>r$ and $b$ a non-zero integer. Then the diophantine equation $bn_1!!\cdots n_r!!=f(x)$ has only finitely many integer solutions.
\end{cor}
\begin{cor}
	Let $b$ be a non-zero integer and $f(x)\in\mathbb{Q}[x]$ a polynomial with factorization 
	\begin{center}
		$f(x)=f_1(x)^{e_1}\cdots f_u(x)^{e_u}$,
	\end{center}
	where the  $f_i(x)$ are irreducible polynomials of degree $d_i$. If $d_i\geq2$ and $d_1e_1+\cdots +d_ue_u>r$ or if $\mathrm{min}\{d_1e_1,...,d_ue_u\}>r$ then the equation $bn_1!!\cdots n_r!!=f(x)$ has only finitely many integer solutions.
\end{cor}
Assuming the ABC-conjecture we show: 
\begin{thm}
	Let $f(x)\in\mathbb{Q}[x]$ be a polynomial of degree $\geq 2$. Then the ABC-conjecture implies that $n!!=f(x)$ has only finitely many integer solutions.	
\end{thm}
\begin{thm}
	Let $f(x)\in\mathbb{Q}[x]$ be a polynomial of degree $\geq 2$ which is not monomial and has at least two distinct roots and let $b$ be a non-zero integer. Then the ABC-conjecture implies that $bn_1!!\cdots n_r!!=f(x)$ has only finitely many integer solutions.	
\end{thm}
\begin{prob}
	\textnormal{Let $A_i$ be fixed non-zero integers. Study the diophantibe equations}
	\begin{center}
		$\sum_{i=1}^rA_in_i!=f(x)$ \textnormal{and} $\sum_{i=1}^rA_in_i!=f(x,y)$.
	\end{center}
\end{prob}
\begin{prob}
	\textnormal{Let $A_i$ be fixed non-zero integers and $S_i\subset \mathbb{Z}$ infinite subsets. Study the diophantibe equations}
	\begin{center}
	$\sum_{i=1}^rA_in_{S_i}!=f(x)$ \textnormal{and} $\sum_{i=1}^rA_in_{S_i}!=f(x,y)$.
\end{center}
\end{prob}
\begin{prob}
	\textnormal{Let $A_i$ be fixed non-zero integers and $S_i\subset \mathbb{Z}$ infinite subsets. Study the diophantibe equations}
	\begin{center}
		$\prod_{i=1}^rA_in_{S_i}!=f(x)$ \textnormal{and} $\prod_{i=1}^rA_in_{S_i}!=f(x,y)$.
	\end{center}
\end{prob}
%\begin{prob}
%	\textnormal{gcshcg}
%\end{prob}
\noindent
{\small \textbf{Acknowledgement}. I am grateful to Wataru Takeda for some helpful explanations.}

\noindent

\section{Proof of Theorems 1.1 and 1.2}
First note that after multiplying the equation $An!+Bm!=f(x,y)$ by a certain integer we can assume $f(x,y)\in\mathbb{Z}[x,y]$. 
To prove the statement of Theorem 1.1 we actually go through the proof of Theorem 4.1 in \cite{WT} and through the proof of \cite{EO}, Satz 1.4 and adapt both to our situation. We follow the notation of \cite{WT}. Let $f(x,y)=a_dx^d+a_{d-1}x^{d-1}y+\cdots + a_{d-1}xy^{1}+a_0y^d$ be an irreducible polynomial and let $K_F$ be the splitting field of $f(x,1)$. Denote by $C_F$ the set of conjugacy classes of the Galois group $G_F=\mathrm{Gal}(K_F/\mathbb{Q})$ whose cycle type $[h_1,...,h_s]$ satisfies $h_i\geq 2$. For a cycle $\sigma$, the cycle type is defined as the ascending ordered list $[h_1,...,h_s]$ of the sizes of cycles in the cycle decomposition of $\sigma$. For further details we refer to Chapters 2, 3 and 4 in \cite{WT}. Of particular interest are the proofs of Lemma 2.1, Lemma 3.1, Theorem 3.6 and Theorem 4.1. We proceed with our proof. Since $d\geq 2$, we conclude from \cite{WT}, Lemma 2.1 that $C_F\neq \emptyset$. Now let $C\in C_F$ be a fixed conjugacy class of $G_F$. We say that a prime $p$ corresponds to $C$ if the Frobenius map $(p,K_F/\mathbb{Q})$ belongs to $C$ (see \cite{WT}, chapter 2 for details). Let $g=\mathrm{gcd}(a_d,...,a_0)$ and $N=gp_1\cdots p_uq_1^{l_1}\cdots q_v^{l_v}$ where $q_i$ are primes corresponding to a conjugacy class in $C_F$ satisfying $\mathrm{gcd}(q,a\Delta_{mod})=1$ where $a\in\{a_d,a_0\}$ and $p_j$ are the other primes (see \cite{WT}, Lemma 3.1 for details). Here $\Delta_{mod}$ denotes a certain modified discriminant of $f(x,1)$ and is defined by 
\begin{center}
	$\Delta_{mod}=\frac{\Delta_{f(x,1)}}{\mathrm{gcd}(a_n,...,a_0)^{2n-2}}$.
\end{center}
The assumption $d\geq 2$ and \cite{WT}, Lemma 3.1 then imply that if $N$ is represented by $f(x,y)$ and $q|N$ for a prime $q$ corresponding to $C$ satisfying $\mathrm{gcd}(q,a\Delta_{mod})=1$, then $N$ is divisible by $q^d$ at least. Assume $m>2\mathrm{max}\{|A|,|B|\}$. Now if $q<m<2q$ and if $n>2m$, then there is no solution to the equation $An!+Bm!=f(x,y)$. Indeed, $q<m<2q$ implies $\frac{m}{2}<q<m$. Therefore, $q$ appears with exponent one in $Bm!$. Now if $n>2m$, the prime $q$ appears in $An!$ with exponent two. This shows that $An!+Bm!$ is divided by $q$ but not by $q^2$. And since $d\geq 2$, we conclude that at least $q^2$ divides $f(x,y)$.  

%Without loss of generality, let $n_r\geq n_{r-1}\geq \cdots \geq n_1$ and assume $n_r$ is big enough, say $n_r>2\mathrm{max}\{|b|,A_1,...A_r\}$. Since the second smallest positive integer divisible by $p$ is $2p$, there is no solution for $bn_1!A_1^{n_1}\cdots n_r!A_r^{n_r}=f(x,y)$ if $p<n_r<2p$. Indeed, since $p<n_r<2p$ we see that $2p<2n_r$ and hence $p<n_r<2p<2n_r$. This implies $\frac{n_r}{2}<p<n_r$. Since  $n_r>2\mathrm{max}\{|b|,A_1,...A_r\}$ and since $d>r$, we see that $p^d$ does not divide $bn_1!A_1^{n_1}\cdots n_r!A_r^{n_r}$.

 Now apply \cite{WT}, Theorem 3.6 (as in the proof of Theorem 4.1 in \emph{loc.cit.}) to conclude that there exists a prime $q'$ corresponding to $C$ with $q'\in (q,2q)$ and satisfying $\mathrm{gcd}(q',a\Delta_{mod})=1$. Therefore, by the same argument as before, if $q'<m<2q'$ and $n>2m$ there are no integer solutions for $An!+Bm!=f(x,y)$. By induction we conclude that whenever $m>2\mathrm{max}\{|A|,|B|\}$ and $m>q$ and $n>2m$ there are no integer solutions. So there can be solutions if $m\leq 2\mathrm{max}\{|A|,|B|\}$ or $m\leq q$ or if $n\leq 2m$. We look at each of the seven possible cases. In all the cases where $m\leq p$ and $ n\leq 2m$ we are done since $n,m$ are bounded. 
 \begin{itemize}
 \item[1)]The case  $m>2\mathrm{max}\{|A|,|B|\}$, $m>q$ and $n\le 2m$: suppose we have infinitely many integer solutions $(n,m,x,y)$ to $An!+ Bm!=f(x,y)$. If there are only finitely many such $n,m$ but infinitely many $x,y$ we are done. So we assume there are infinitely many $m$. In this case we can argue as in the proof of Satz 4 in \cite{EO}. One can use  \cite{WT}, Theorem 3.6 to conclude that\\
 \textbf{claim}:\\
 for large $m$ there is allways a prime $p$ corresponding to $C$ and satisfying $\mathrm{gcd}(p,a\Delta_{mod})=1$ such that 
 \begin{center}
 	$\frac{m}{2}<p<\frac{m}{2}+\frac{m}{12\mathrm{log}m}$.
 \end{center}
\textbf{proof of claim}:\\
as above, \cite{WT}, Theorem 3.6 implies that there is a prime $q'$ corresponding to $C$ with $q'\in (q,2q)$ and satisfying $\mathrm{gcd}(q',a\Delta_{mod})=1$. By induction, we can assume that there is a prime $q_0$ corresponding to $C$ and satisfying $\mathrm{gcd}(q_0,a\Delta_{mod})=1$ such that $m/2\geq q_0$. Applying \cite{WT}, Lemma 3.5 and Theorem 3.6 again, we conclude that for any $D>0$ there is a prime $p$ corresponding to $C$ with $p\in (q_0,Dq_0)$ and satisfying $\mathrm{gcd}(p,a\Delta_{mod})=1$. As long as $q_0\leq m/2<p$, we see that 
\begin{center}
	$\frac{m}{2}<p<\frac{Dm}{2}<\frac{m}{2}+\frac{m}{12\mathrm{log}m}$.
\end{center}
Note that $D>0$ can be chosen small enough such that the last inequality holds. If $m/2\geq p$, we use the same argument to produce a a prime $p'$ corresponding to $C$ with $p'\in (p,Dp)$ and satisfying $\mathrm{gcd}(p',a\Delta_{mod})=1$. Again, as long as $p\leq m/2<p'$, we have 
\begin{center}
	$\frac{m}{2}<p'<\frac{Dm}{2}<\frac{m}{2}+\frac{m}{12\mathrm{log}m}$.
\end{center}
Now by an induction argument, we get the claim.\\

With the help of the claim, we see that 
 \begin{center}
 $m<2p<m+\frac{m}{6\mathrm{log}m}$.
 \end{center}
As above, if $n>m+\frac{m}{6\mathrm{log}m}$ we can coclude that $p$ divides $An!+Bm!$ but not $p^2$. Since $p^2$ divides $f(x,y)$, there can be no solutions in this case. Therefore, we assume  
 \begin{center}
 	$n\leq m+\frac{m}{6\cdot \mathrm{log}m}$.
 \end{center}
Since we assumed $n>m$, we can consider
\begin{center}
	$Bm!\bigl(\frac{An!}{Bm!} + 1 \bigr)=f(x,y)$.
\end{center}
Then $Bm!$ contains any prime $p$ with $\frac{m}{2}<p\leq m$ with exponent exactly one. Therefore, $\frac{An!}{Bm!} + 1$ must contain these primes, too. Using the prime number theorem one can show that for $m$ big enough the product of these primes is bigger than $2^{m/2}+1$. Hence 
\begin{center}
	$|\frac{An!}{Bm!}|>2^{m/2}$.
\end{center}
From $n\leq m+\frac{m}{6\cdot \mathrm{log}m}<2m$, it follows
\begin{center}
	$|\frac{An!}{Bm!}|<\frac{|A|}{|B|}n^{(n-m)}<\frac{|A|}{|B|}(2m)^{\frac{m}{6\cdot \mathrm{log}m}}<\frac{|A|}{|B|}2^{\frac{m}{6\cdot \mathrm{log}m}}e^{\frac{m}{6}}$.
\end{center}
But for large $m$ this contradicts 
\begin{center}
	$|\frac{An!}{Bm!}|>2^{m/2}$.
\end{center}
This implies that here can be only finitely many $m$ and hence finitely many $n$.
 \item[2)] The case $m>2\mathrm{max}\{|A|,|B|\}$, $m\leq q$ and $n>2m$: Assume $n\geq 2q$ otherwise $n$ and $m$ are bounded. We see that the exponent of $q$ in $Bm!$ is at most one, whereas the exponent of $q$ in $An!$ is at least two. The exponent of $q$ in $f(x,y)$ is also at least two. Hence there are no solutions in this case. 
 
 \item[3)] The case  $m\leq 2\mathrm{max}\{|A|,|B|\}$, $m\leq q$ and $n>2m$: As above, we can find a prime $q'$ corresponding to $C$ with $q'\in (q,2q)$ and satisfying $\mathrm{gcd}(q',a\Delta_{mod})=1$. We can repeat that argument until we find a prime $q_0>2\mathrm{max}\{|A|,|B|\}$ corresponding to $C$, satisfying $\mathrm{gcd}(q_0,a\Delta_{mod})=1$. Now we assume $n>2q_0$ otherwise $n$ and $m$ are bounded. Now since $m<q_0$ and $q_0>2\mathrm{max}\{|A|,|B|\}$ and $n>2q_0$, we see that the exponent of $q_0$ in $Bm!$ is zero whereas the exponent in $An!$ and $f(x,y)$ is at least two. Consequently, there are no solutions in this case.
 \item[4)] The case  $m\leq 2\mathrm{max}\{|A|,|B|\}$, $m>q$ and $n>2m$: we argue as in case 3). Thus we can find a prime $q_0$ corresponding to $C$ and satisfying $\mathrm{gcd}(q_0,a\Delta_{mod})=1$ such that $m\leq q_0$. Then assume $n>2q_0$ and continue as in 3).

 %The case $m\leq p$: relevant is just the case $m\leq p$ where $n>2m$. Then either $n\leq p$ or $n>p$. If $n\leq p$ we are done. So let us assume $n>p$. If $p<n<2p$, then $n$ is bounded and we are done. So let us assume $n\geq 2p$. Then the exponent of $p$ in $n!$ is at least two.. On the other hand $m\leq p$ implies that the exponent of $p$ in $m!$ is at most one. Then $n!\pm m!$ is divisible by $p$ but not by $p^2$. On the other hand $p$ devides $f(x,y)$ with exponent $d\geq 2$. Therefore, there are no integer solutions in this case.
 \end{itemize}
Summarizing we see that there are only finitely many $n,m$ such that $An!+ Bm!$ is represented by $f(x,y)$. Now we use Thue's theorem to conclude that $An!+Bm!=f(x,y)$ has finitely many integer solutions if the degree of $f(x,y)$ is $\geq 3$.\\
%This shows that there are only finitely many $(n_1,...,n_r)$ such that $bn_1!A_1^{n_1}\cdots n_r!A_r^{n_r}$ is represented by $f(x,y)$. Now if $r\geq 2$, we have $d\geq 3$ and we can use Thue's theorem to conclude that there are indeed only finitely many integer solutions.\\
Proof of Theorem 1.2:\\
If all $d_i\geq 2$, then proceed as follows: let $K_{F_j}$ be the splitting field of $f_j(x,1)$ and denote by $C_{F_j}$ be the set of conjugacy classes whose cycle type has sizes $\geq 2$. Now let $q$ be a prime corresponding to a conjugacy class $C\in \cap_{j=1}^uC_{F_j}$ satisfying $\mathrm{gcd}(q,a\Delta_{mod})=1$. Assume $N$ is represented by $f(x,y)$. Then there are $x_0,y_0$ such that $q|f(x_0,y_0)$. Hence there exists a polynomial $f_j(x,y)=a_{j,d_j}x^{d_j}+\cdots + a_{j,0}y^{d_j}$ such that $q|f_j(x_0,y_0)$. Now $\mathrm{gcd}(q,a\Delta_{mod})=1$ with  $a\in\{a_d,a_0\}$ implies $\mathrm{gcd}(q,a(j)\Delta_{j, mod})=1$ where $a(j)\in\{a_{j,d_j},a_{j,0}\}$. Here $\Delta_{j, mod}$ denotes the modified discriminant of $f_j(x,1)$. It follows from \cite{WT}, Lemma 3.1 that $q^d|f(x,y)$. The rest of the proof is the same as in the proof of Theorem 1.1. Now let $\mathrm{min}\{d_1e_1,...,d_ue_u\}=d_{i_0}e_{i_0}$ and assume $d_{i_0}=1$ and $e_{i_0}>1$. Then one can argue as in \cite{WT}, Theorem 4.7 to conclude that, let's say, $q$ divides $f(x,y)$ at least $e_{i_0}>1$ times. Again, the rest of the proof is as in the proof of Theorem 1.1 above.
%The proof of Theorem 1.2 follows exactly the lines of the proofs of Theorems 4.3 and 4.7 in \cite{WT} and the proof of Theorem 1.1 from above.% As in the proof of Theorem 1.1 above, we assume without loss of generality that $n_r\geq \cdots \geq n_1$ and choose $n_r$ big enough. 
%Again, one concludes the statement from an inductions argument as in the proof of Theorem 1.1. 

%\section{Proof of Corollary 1.4}

\section{Proof of Theorems 1.5 and 1.6}
For a better readibility, to keep the arguments clear and to avoid a lot of indices, we formulate the proof only for $r=3$. First note that after multiply the equation $bn_1!!\cdots n_r!!=f(x,y)$ with a certain integer we can assume that $f(x,y)\in\mathbb{Z}[x,y]$. Secondly, note that if $n=2m$ is even, we have $n!!=2^mm!$. So in the case $n_,n_2,n_3$ are all even, we conclude form \cite{NO}, Theorem 1.5 that there are finitely many $(n_1,n_2,n_3)$ such that $bn_1!!n_2!!n_3!!$ is represented by $f(x,y)$. Now consider the case $n_1,n_2,n_3$ are all odd. Note that for $n_i=2m_i+1$ one has
\begin{center}
	$(2m_i+1)!!=\frac{(2m_i+1)!}{2^{m_i}m_i!}$
\end{center}
We rewrite our diophantine equation and consider
\begin{center}
	$b(2m_1+1)!(2m_2+1)!(2m_3+1)!=2^{m_1}m_1!2^{m_2}m_2!2^{m_3}m_3!f(x,y)$.
\end{center}
As in the proof of Theorem 1.1, let $K_F$ be the splitting field of $f(x,1)$ and denote by $C_F$ the set of conjugacy classes of the Galois group $G_F=\mathrm{Gal}(K_F/\mathbb{Q})$ whose cycle type $[h_1,...,h_s]$ satisfies $h_i\geq 2$. %For a cycle $\sigma$, the cycle type is defined as the ascending ordered list $[h_1,...,h_s]$ of the sizes of cycles in the cycle decomposition of $\sigma$. For further details we refer to Chapters 2,3 and 4 in \cite{WT}. Of particular interest are the proofs of Lemma 2.1, Lemma 3.1, Theorem 3.6 and Theorem 4.1. We proceed with our proof. 
Since $d\geq 2$, we conclude from \cite{WT}, Lemma 2.1 that $C_F\neq \emptyset$. Now let $C\in C_F$ be a fixed conjugacy class of $G_F$. %We say that a prime $p$ corresponds to $C$ if the Frobenius map $(p,K_F/\mathbb{Q})$ belongs to $C$ (see \cite{WT}, chapter 2 for details). 
Again, let $g=\mathrm{gcd}(a_d,...,a_0)$ and $N=gp_1\cdots p_uq_1^{l_1}\cdots q_v^{l_v}$ where $q_i$ are primes corresponding to a conjugacy class in $C_F$ satisfying $\mathrm{gcd}(q,a\Delta_{mod})=1$ where $a\in\{a_d,a_0\}$ and $p_j$ are the other primes (see \cite{WT}, Lemma 3.1 for details). The assumption $d>r>2$ and \cite{WT}, Lemma 3.1 then imply that if $N$ is represented by $f(x,y)$ and $q|N$ for a prime $q$ corresponding to $C$ satisfying $\mathrm{gcd}(q,a\Delta_{mod})=1$, then $N$ is divisible by $q^d$ at least. %Consequently, $f(x,y)=q!$ has no solution. 
Let us assume $m_1\geq m_2\geq m_3$ and $m_1>2|b|$. Now if $q<2m_1+1<2q$, then there is no solution to the equation 
\begin{center}
	$b(2m_1+1)!(2m_2+1)!(2m_3+1)!=2^{m_1}m_1!2^{m_2}m_2!2^{m_3}m_3!f(x,y)$.
\end{center}
Indeed, $q<2m_1+1<2q$ implies $\frac{2m_1+1}{2}<q<2m_1+1$. Therefore, $q$ appears with exponent exactly one in $(2m_1+1)!$. We can assume $m_1>2$. Since $m_3\leq m_2\leq m_1<q$ we see that $q$ does not divide $2^{m_i}m_i!$. Therefore, $q$ must divide $f(x,y)$. As mentioned above, this implies that $q^d$ divides $f(x,y)$. Since $d>r$ by assumption, we conclude that there are no solutions if $q<2m_1+1<2q$. Now by the same induction argument as in the proof of Theorem 1.1 we find that there are no solutions to the above equation whenever $2m_1+1>q$. Now let $n_1,n_2$ be odd and $n_3$ even. Then the diophantine equation becomes
\begin{center}
	$b(2m_1+1)!(2m_2+1)!2^{m_3}m_3!=2^{m_1}m_1!2^{m_2}m_2!f(x,y)$.
\end{center}
Because we have symmetry with respect to $m_1$ and $m_2$, we can consider three cases: assume $m_1\geq m_2\geq m_3$. Then one can apply the arguments from before. The same can be done if $m_1\geq m_3\geq m_2$. The last case is $m_3\geq m_1\geq m_2$. Here we assume $m_3\geq 2|b|$. Now if $q<2m_3+1<2q$ we argue as before. The case where two of the numbers $n_1,n_2,n_3$ are even and one is odd is similar and left to the reader. \\
Proof of Theorem 1.6:\\ 
If all $d_i\geq 2$, then proceed as follows: let $K_{F_j}$ be the splitting field of $f_j(x,1)$ and denote by $C_{F_j}$ be the set of conjugacy classes whose cycle type has sizes $\geq 2$. Now let $q$ be a prime corresponding to a conjugacy class $C\in \cap_{j=1}^uC_{F_j}$ satisfying $\mathrm{gcd}(q,a\Delta_{mod})=1$. Assume $N$ is represented by $f(x,y)$. Therefore, there are $x_0,y_0$ such that $q|f(x_0,y_0)$. Then there exists a polynomial $f_j(x,y)=a_{j,d_j}x^{d_j}+\cdots + a_{j,0}y^{d_j}$ such that $q|f_j(x_0,y_0)$. Now $\mathrm{gcd}(q,a\Delta_{mod})=1$ with  $a\in\{a_d,a_0\}$ implies $\mathrm{gcd}(q,a(j)\Delta_{j, mod})=1$ where $a(j)\in\{a_{j,d_j},a_{j,0}\}$. Here $\Delta_{j, mod}$ denotes the modified discriminant of $f_j(x,1)$. It follows from \cite{WT}, Lemma 3.1 that $q^d|f(x,y)$. The rest of the proof is the same as in the proof of Theorem 1.5. Now let $\mathrm{min}\{d_1e_1,...,d_ue_u\}=d_{i_0}e_{i_0}$ and assume $d_{i_0}=1$ and $e_{i_0}>1$. Then one can argue as in \cite{WT}, Theorem 4.7 to conclude that, let's say, $q$ divides $f(x,y)$ at least $e_{i_0}>1$ times. Again, the rest of the proof is as in the proof of Theorem 1.5 above.

%\section{Proof of Corollary 1.8}

\section{Proof of Theorems 1.9 and 1.10}
We can assume $f(x)\in\mathbb{Z}[x]$. First, if $n=2m$ is even we have $n!!=2^mm!$ and the result follows from \cite{NO}, Theorems 1.1 and 1.2. So let $n=2m+1$. Notice that
\begin{center}
	$(2m+1)!!=\frac{(2m+1)!}{2^mm!}$.
\end{center} 
Let us consider a polynomial 
\begin{eqnarray*}
	f(x)=a_0x^d+a_1x^{d-1}+\cdots +a_d
\end{eqnarray*}
with $a_i\in \mathbb{Z}$. Now multiply the equation $\frac{(2m+1)!}{2^mm!}=f(x)$ by $d^da_0^{d-1}$. We obtain
\begin{eqnarray*}
	y^d+b_1y^{d-1}+\cdots +b_d=c(\frac{(2m+1)!}{2^mm!})
\end{eqnarray*}
for a constant $c$, where $c=bd^da_0^{d-1}$ and $y:=a_0dx$. Notice that $b_i=d^ia_ia_0^{i-1}$ so that we can make the change of variable $z:=y+\frac{b_1}{d}$. %Since we are assuming that $f(x)$ has at least two distinct roots, the change of variable produces a polynomial that does not have a monomial of degree $d-1$. 
Therefore we get the following equation
\begin{eqnarray}
	z^d+c_2d^{d-2}+\cdots +c_d=c(\frac{(2m+1)!}{2^mm!}).
\end{eqnarray}
Notice that $c_i$ are integer coefficients wich can be computed in terms of $a_i$ and $d$. Now let $Q(X)=X^d+c_2X^{d-2}+\cdots +c_d$ and notice that when $|z|$ is large one has
\begin{eqnarray}
	\frac{|z|^d}{2}<|Q(z)|<2|z|^d.
\end{eqnarray}
For the rest of the proof we denote by $C_1,C_2,...$ computable positive constants depending on the coefficients $a_i$ and eventually on some small $\epsilon >0$ which comes into play later by applying the ABC-conjecture.

Whenever $(m,z)$ is a solution to $\frac{(2m+1)!}{2^mm!}=f(x)$ we conclude from (3) and (4) that there exist constants $C_1$ and $C_2$ such that
\begin{eqnarray}
	|d\cdot\mathrm{log}|z|-\mathrm{log}(\frac{(2m+1)!}{2^mm!})|<C_1,
\end{eqnarray}
for $|z|>C_2$ (see \cite{L} equation (10)). %For technical reasons, we assume that $C_2$ is large enough with respect to $C_1$. 
Now let $R(X)\in \mathbb{Z}[X]$ be such that $Q(X)=X^d+R(X)$. %By assumption, we may assume that $R(X)$ is non-zero.
If $R(X)=0$, the diophantine equation becomes $c(\frac{(2m+1)!}{2^mm!})=x^d$ which is equivalent to $c(2m+1)!=2^mm!x^d$. Let us assume $m>2|c]$. Then there is a prime $p$ in the interval $((2m+1)/2,2m+1)$ that divides $(2m+1)!$ with exponent one. Since $p>m$ and since $m>2|c|$ we see that $p$ does not divide $2^mm!$. Therefore the exponent of $p$ in $2^mm!x^d$ is at least two. Hence there are no solutions if $m>2|c|$.

%Since $n>2m$, the number $p^2$ devides $n!$. Therefore $c(n!\pm m!)$ is devided by $p$ but not by $p^2$. Since $d\geq 2$ there are no solutions. Now we consider the case $m>2|c|$ and $n\leq 2m$. If there are infinitely many solutions to $c(n!\pm m!)=x^d$ we could find solutions with arbitrary big $m$. The prime number theorem yields 
%\begin{center}
%	$n\leq m+\frac{m}{6\cdot \mathrm{log}m}$.
%\end{center}
%From
%\begin{center}
%	$cm!\bigl(\frac{n!}{m!}\pm 1\bigl)=x^d$
%\end{center}
%we conclude that $cm!$ contains a prime  $q$ with $m/2<q\leq m$ with exponent exactly one. Now apply the same arguments as in \cite{EO},p.255 to construct a contradiction. The last case to consider is $m\leq 2|c|$. Then either $n\leq 2m$ or $n>2m$. If $n\leq 2m$ we are done, since there are obviously only finitely many solutions. So let $n>2m$. Since we only have finitely many $m$ in the set of all solutions we can consider the equation $c(n!\pm m!)=x^d$ and notice that the equation is actually $cn!=x^d\pm cm!$. For any of the finitely many $m$ we therefore have a diophantine equation of the form $cn!=P(x)$ where $P(x)$ is a polynomial of degree $\geq 2$. Now Luca's result \cite{L} states that ABC implies that there are only finitely many integer solutions. 

So let us assume that $R(X)$ is non-zero and let $j\leq d$ be the largest integer with $c_j\neq 0$. We rewrite (3) as
\begin{eqnarray*}
	z^j+c_2z^{j-2}+\cdots +c_j=\frac{c(\frac{(2m+1)!}{2^mm!})}{z^{d-j}}.
\end{eqnarray*}
Let $R_1(X)$ be the polynomial 
\begin{eqnarray*}
	R_1(X):= \frac{R(X)}{X^{d-j}}=c_2X^{j-2}+\cdots +c_j.
\end{eqnarray*}
It is shown in \cite{L} there are constants $C_3$ and $C_4\geq C_2$ such that
\begin{eqnarray*}
	0<|R_1(z)|< C_3|z|^{j-2},
\end{eqnarray*}
for $|z|> C_4$. So let us write 
\begin{center}
	$z^j+R_1(z)=\frac{c(\frac{(2m+1)!}{2^mm!})}{z^{d-j}}$. 
\end{center}
For $D=\mathrm{gcd}(z^j, R_1(z))$ we get
\begin{eqnarray*}
\frac{z^j}{D}+\frac{R_1(z)}{D}=\frac{c(\frac{(2m+1)!}{2^mm!})}{Dz^{d-j}}
\end{eqnarray*}
Applying the ABC-conjecture to $A=\frac{z^j}{D}$, $B=\frac{R_1(z)}{D}$ and $C=\frac{c(2m+1)!}{z^{d-j}D2^mm!}$, we find
\begin{eqnarray}
	\frac{|z|^j}{D}< C_5N(\frac{z^jR_1(z)c(2m+1)!}{D^3z^{d-j}2^mm!})^{1+\epsilon},
\end{eqnarray}
where $C_5$ depends only on $\epsilon$. It is shown in \cite{L}, p.272 that 
\begin{eqnarray}
	N(\frac{|z|^j}{D})\leq |z|,\\
	N(\frac{R_1(z)}{D})<\frac{C_3|z|^{j-2}}{D}.
\end{eqnarray}
Moreover, we have
\begin{eqnarray*}
	N(\frac{c(2m+1)!}{z^{d-j}D2^mm!})\leq N(c)N((2m+1)!).
\end{eqnarray*}
%Now there is a $M\in \mathbb{N}$ such that $A_i^{n_i}\leq n_i!$ for all $n_i>M$. 
%Lemma 3.1 yields 
%\begin{eqnarray*}
%	N(\frac{c(n!\pm m!)}{z^{d-j}D})\leq N(c)N(n!\pm m!)=N(c)N(m!(\frac{n!}{m!}\pm 1))\leq C_6N(m!)\cdot N((2n)!),
%\end{eqnarray*} 
%where $C_6=N(c)$. 
From (1) it follows
\begin{eqnarray}
	N(\frac{c(2m+1)!}{z^{d-j}D2^mm!})<C_64^{2m+1}
\end{eqnarray}
and from (7), (8) and (9) we get
\begin{eqnarray}
	N(\frac{|z|^j}{D})N(\frac{R_1(z)}{D})N(\frac{c(2m+1)!}{z^{d-j}D2^mm!})<\frac{C_3C_6|z|^{j-1}4^{(2m+1)}}{D}
\end{eqnarray}
From inequalities (6) and (10), we obtain
\begin{eqnarray}
	\frac{|z|^j}{D}<C_7\bigl(\frac{|z|^{j-1}4^{(2m+1)}}{D}\bigr)^{(1+\epsilon)}
\end{eqnarray}

If we choose $\epsilon =\frac{1}{2d}\leq \frac{1}{2j}$, inequality (14) implies that
\begin{eqnarray*}
	|z|^{1/2}<|z|^{1+\epsilon -\epsilon j}< C_84^{(2m+1)(1+\epsilon)},
\end{eqnarray*}
or simply 
\begin{eqnarray*}
	\mathrm{log}|z|<C_9(2m+1)+C_{10}.
\end{eqnarray*}
Thus
\begin{eqnarray*}
	d\cdot\mathrm{log}|z|<C_{11}m+C_{12}.
\end{eqnarray*}
This gives
\begin{eqnarray}
	\mathrm{log}|(\frac{(2m+1)!}{2^mm!})|<C_1+d\cdot \mathrm{log}|z|< C_{11}m+C_{13}.
\end{eqnarray}
%We can simplify (15) and finally obtain
%\begin{center}
%$\mathrm{log}(n_1!)+\mathrm{log}(n_2!)+\cdots +\mathrm{log}(n_r!)<C_{14}n_1+C_{14}n_2+\cdots C_{14}n_r+C_{13}$.	
%	\end{center}
Now we can conclude that only finitely many $m$ satisfy (12). For the convenience to the reader we give an argument. So let us consider an inequality of the form
\begin{center}
	$\mathrm{log}|\frac{(2m+1)!}{2^mm!}|<A'm+C'$
\end{center}
where $A'$ and $C'$ are positive constant intergers. We first rewrite the equation by applying the $e$-function. This gives 
\begin{center}
	$\frac{(2m+1)!}{2^mm!}<E\cdot e^{A'm}$.
\end{center}
The Stirling formula then yields
\begin{center}
	$\frac{(\frac{2m+1}{e})^{2m+1}}{2^mm^m}<\frac{(\frac{2m+1}{e})^{2m+1}}{2^mm!}<E\cdot e^{A'm}$.
\end{center}
And since
\begin{center}
	$\frac{(\frac{2m+1}{e})^{2m+1}}{2^mm^m}=\frac{(2m+1)^{2m+1}}{2^me^{2m+1}m^m}$
\end{center}
we get the inequality
\begin{center}
$\frac{(2m+1)^{2m+1}}{m^m}<E\cdot e^{A'm}e^{2m+1}2^m=E'\cdot e^{(A'+2+ln(2))m}=E'\cdot e^{B'm}$
\end{center}
%Assume the equation $n!+m!=f(x)$ has infinitely many solutions. Then we must have infinitely many pairs $(n,m)$. There are three cases where these pairs can lie in the lattice of natural points. 
%Assume there are infinitely many pairs $(n,m)$ of natural numbers satisfying the inequality. There are three cases: 
%\begin{itemize}
%	\item[1)] there are infinitely many pairs with $n=m$. Then the above inequality simplyfies to 
%	\begin{center}
%		$2(\frac{n}{e})^n<E\cdot e^{(A+B)n}$.
%	\end{center} 
%	But this inequality is satisfied only for finitely many $n$. This gives a contradiction.
%	\item[2)] there are infinitely many pairs $(n,m)$ with $n<m$. In this case the inequality becomes
%	\begin{center}
%		$(\frac{n}{e})^n+(\frac{m}{e})^m<E\cdot e^{(An+Bm)}<E\cdot e^{(A+B)m}$.
%	\end{center}
%	But this inequality is valid only for finitely many $m$. This gives a contradiction.
	
%	\item[3)] there are infinitely many pairs $(n,m)$ with $m<n$. The same argument as in 2) gives a contradition. 
%\end{itemize}
%The argument for the $|n!-m!|<An+Bm+C$ is similar. 
But this inequality is satified only for finitely many $m$. We finally conclude from (5) that $|z|< C_{15}$. This completes the proof of Theorem 1.9.\\
Proof of Theorem 1.10:\\
We formulate the proof only for $r=2$ since the proof for arbitrary $r>2$ is analogous. So let us consider the equation $bn!!m!!=f(x)$. Note that if $d=2$ and $b\neq 0$ there could be infinitely many solutions. This happens for instance if $f(x)=x^2$ and follows from \cite{NO}, Theorem 1.1. Assuming that $f(x)$ is not monomial and has at least two distinct roots enables us to make a change of variables as in the proof of Theorem 1.9. Therefore, we can assume $f(z)=z^d+c_2d^{d-2}+\cdots +c_d$, with $c_2d^{d-2}+\cdots +c_d\neq 0$. If $n=2n_1,m=2m_1$ are both even, we actually have the equation
\begin{center}
	$c2^{n_1}n_1!2^{m_1}m_1!=f(z)$
\end{center}
for some integer constant $c$ and the result follows from \cite{NO}, Theorem 1.2. Now let us consider the case where $n=2n_1+1$ and $m=2m_1+1$ are both odd. Then we have the equation
\begin{center}
	$f(z)=c\frac{(2n_1+1)!}{2^{n_1}n_1!}\frac{(2m_1+1)!}{2^{m_1}m_1!}$.
\end{center}
Let $R_1(X), A, B, C$ and $D$ be as in the proof of Theorem 1.9. For the algebraic radical we then find
\begin{center}
	$N(\frac{c}{z^{d-j}D}\frac{(2n_1+1)!}{2^{n_1}n_1!}\frac{(2m_1+1)!}{2^{m_1}m_1!})<\tilde{C}4^{2(n_1+m_1)+2}$
\end{center}
where $\tilde{C}$ is a certain positive constant integer. As in the proof of Theorem 1.9, we finally get 
\begin{center}
	$\mathrm{log}(|\frac{(2n_1+1)!}{2^{n_1}n_1!}\frac{(2m_1+1)!}{2^{m_1}m_1!}|)<A^*n_1+B^*m_1+C^*$
\end{center}
where $A^*, B^*$ and $C^*$ are positive constant integers. Applying $e$-function and Stirling formula finally yields
\begin{center}
	$\frac{(2n_1+1)^{2n_1+1}}{n_1^{n_1}}\frac{(2m_1+1)^{2m_1+1}}{m_1^{m_1}}<E^*e^{A^*n_1}e^{B^*m_1}$.
\end{center}
This inequality however is satisfied only for finitely many $(n_1,m_1)$. The case where $n$ is odd and $m$ is even is analogous. We omit the proof. 

\vspace{0.5cm}
\noindent
{\tiny HOCHSCHULE FRESENIUS UNIVERSITY OF APPLIED SCIENCES 40476 D\"USSELDORF, GERMANY.}
E-mail adress: sasa.novakovic@hs-fresenius.de\\
\noindent
{\tiny MATHEMATISCHES INSTITUT, HEINRICH--HEINE--UNIVERSIT\"AT 40225 D\"USSELDORF, GERMANY.}
E-mail adress: novakovic@math.uni-duesseldorf.de

\end{document}